\documentclass[12pt]{amsart}

\usepackage{amssymb}
\usepackage{bbm}
\usepackage{mathtools}
\usepackage{enumitem}
\usepackage{tikz-cd}
\usepackage{calc}
\usepackage{adjustbox}
\usepackage{hyperref}

\makeatletter
\renewcommand*\env@matrix[1][\arraystretch]{%
\edef\arraystretch{#1}%
\hskip -\arraycolsep
\let\@ifnextchar\new@ifnextchar
\array{*\c@MaxMatrixCols c}}
\makeatother

\newtheoremstyle{break}
  {\topsep}{\topsep}
  {\itshape}{}
  {\bfseries}{}
  {\newline}{}

\theoremstyle{plain}
\newtheorem{thm}{Theorem}[section]
\newtheorem{lem}[thm]{Lemma}
\newtheorem{prop}[thm]{Proposition}
\newtheorem{cor}{Corollary}

\theoremstyle{definition}
\newtheorem{defn}{Definition}[section]
\newtheorem{conj}{Conjecture}[section]
\newtheorem{exmp}{Example}[section]

\theoremstyle{remark}
\newtheorem{rem}{Remark}
\newtheorem{note}{Note}

\usepackage{framed}
{\end{quotation}\end{leftbar}}

\begin{document}

\title[On the representation theory of non-semisimple (graded) deformed Fomin-Kirillov algebras.]{On the representation theory of non-semisimple (graded) deformed Fomin-Kirillov algebras.}
\author{A. Alia}
\author{I. Heckenberger}
\address{A. Alia: Philipps-Universit\"at Marburg, FB Mathematik und Informatik, Hans-Meer\-wein-Stra\ss e, 35032 Marburg, Germany.}
\email{alia@mathematik.uni-marburg.de}
\address{I. Heckenberger: Philipps-Universit\"at Marburg, FB Mathematik und Informatik, Hans-Meer\-wein-Stra\ss e, 35032 Marburg, Germany.}
\email{heckenberger@mathematik.uni-marburg.de}

\begin{abstract}
This work is motivated to study the representation theory of the non-semisimple deformed Fomin-Kirillov algebras $\mathcal{D}_4(\alpha_1, \alpha_2)$. In particular, we consider Gabriel's theorem applications in regard of constructing algebraic presentations.
\end{abstract}

\maketitle
\tableofcontents

\section*{Introduction}

Here and throughout this paper we assume $K$ to be an algebraically closed field of zero characteristic, we denote by $\mathbb{S}_n$ the symmetric group on $n$ letters. \\ 

In the context of studying Schubert calculus, Fomin and Kirillov introduced a family of quadratic $K$-algebras $\mathcal{E}_n$, that contains a commutative subalgebra isomorphic to the cohomology ring of the flag manifold. Commonly known as Fomin-Kirillov algebras, we present:

\begin{defn} [\cite{00} Definition 2.1]
Given a positive integer $n \geq 3$. $\mathcal{E}_n$ is the quadratic $K$-algebra generated by $x_{ij} = - x_{ji}$ for $1 \leq i \neq j \leq n$ subject the following relations: 
\begin{subequations}
\begin{align}
x^2_{ij} &= 0 && \mid 1 \leq i, j \leq n \text{ distinct}, \label{eq.00.01} \\
x_{ij}x_{kl} -x_{kl}x_{ij}&= 0 && \mid 1 \leq i , j , k , l \leq n \text{ distinct}, \label{eq.00.02}\\
x_{ij}x_{jk} + x_{jk}x_{ki} + x_{ki}x_{ij} &= 0 && \mid 1 \leq i , j , k \leq n \text{ distinct}. \label{eq.00.03}
\end{align}
\end{subequations}
\end{defn}

While on the surface, it is rather straightforward to present $\mathcal{E}_n$, some of the structure's elementary properties remain challenging to approach: 

\begin{quote}
For example, it is well-known that $\mathcal{E}_n$ is finite-dimensional for $n \leq 5$ and the opposite is \textit{conjectured} to true otherwise.
\end{quote}

Moreover, the nature of $\mathcal{E}_n$ as a braided Hopf algebra over the symmetric group $\mathbb{S}_n$ indicates a strong connection to Nichols-algebras over braided vector spaces:

\begin{quote}
For $n \leq 4$ in \cite{01}, and $n= 5$ in \cite{02} it was showed that the $\mathcal{E}_n$ is a Nichols-algebra. A statement that is conjectured to be true for $n \geq 6$ in \cite{01}, \cite{04}.
\end{quote}

It is also of worth to mention that the algebra $\mathcal{E}_n$ appears to share some distinctive properties with other types of algebras, most famous of which is that of preprojective type of $A_{n-1}$, which shares the same number of indecomposable modules with $\mathcal{E}_n$ for $n \leq 5$ and is known to be of infinite representation-type otherwise; This is recognized as Majid’s conjecture, which does not have a precise expression, nonetheless highlights that the numerology is not accidental, further details can be explored in \cite{05} and in \cite{06}. \\

With that in mind, we remark that from the point of view of graded algebras, $\mathcal{E}_n$ remains an approachable candidate of an algebraic structure that is both naturally and symmetrically graded, moreover, the unique action of the symmetric group $\mathbb{S}_n$ on $\mathcal{E}_n$ defined as:
\begin{align*}
\sigma(xy) = (\sigma x)(\sigma y) : x, y \in \mathcal{E}_n, \sigma \in \mathbb{S}_n &\; \sigma(x_{ij}) = x_{\sigma(i) \sigma(j)} &\mid i \neq j.
\end{align*} 
validates- above other reasons- the study of $\mathcal{E}_n$ from the viewpoint of PBW-deformations. We recall that a PBW deformation of a graded algebra $A$ is a filtered algebra $D$ such that the associated graded algebra of $D$ is isomorphic to $A$, we formalize a definition:

\begin{defn}
Given $\alpha_1, \alpha_2 \in K$. The \textit{deformed Fomin-Kirillov} algebra, denoted by $\mathcal{D}_n(\alpha_1, \alpha_2)$, is the quadratic $K$-algebra generated by $x_{ij} = - x_{ji}$ for $1 \leq i \neq j \leq n$ subject the following relations:
\begin{subequations}
\begin{align}
x^2_{ij} &= \alpha_1 && \mid 1 \leq i, j \leq n \text{ distinct},  \label{eq.00.04} \\
x_{ij}x_{kl}-  x_{kl}x_{ij}&= 0 && \mid 1 \leq i , j , k , l \leq n \text{ distinct},  \label{eq.00.05} \\
x_{ij}x_{jk} + x_{jk}x_{ki} + x_{ki}x_{ij} & = \alpha_2 && \mid 1 \leq i , j , k \leq n \text{ distinct}. \label{eq.00.06}
\end{align}
\end{subequations}
\end{defn}

In 2018, motivated by understanding Nichols and Fomin-Kirillov algebras by means of PBW-deformations. Heckenberger and Vendramin established a framework objected to the classification and the study of representation theory of non-semisimple deformations of Fomin-Kirillov algebras. In particular, the authors theorized:

\begin{thm}\label{thm.00.01} [\cite{03} Theorem 2.11]
The algebra $\mathcal{D}_3(\alpha_1, \alpha_2)$ is semi-simple if and only if:
\begin{align*}
(3\alpha_1-\alpha_2)(\alpha_1+\alpha_2) \neq 0
\end{align*}
In this case $\mathcal{D}_3(\alpha_1, \alpha_2) \cong (K^{2.2})^3$
\end{thm}

Further, proposed:

\begin{prop}\label{thm.00.02}[\cite{03} Proposition 2.15, 2.16]
The following hold:
\begin{enumerate}
\item The algebra $\mathcal{D}_3(\alpha_1, -\alpha_1)$ is isomorphic to the product of three copies of the preprojective algebra of type $A_{2}$. 
\item The algebra $\mathcal{D}_3(\alpha_1, 3\alpha_1)$ is isomorphic to the path algebra of the double Kronecker quiver bounded by the relations of the coinvariant ring of $\mathbb{S}_3$
\end{enumerate}
\end{prop}

Later that year, Wolf in \cite{07} continued the study by examining the case of $\mathcal{D}_4(\alpha_1, \alpha_2)$, where it was proved that the algebra $\mathcal{D}_4(\alpha_1, \alpha_2)$ is semisimple if:
\begin{align*}
\alpha_1(3\alpha_1-\alpha_2)(\alpha_1-\alpha_2)(\alpha_1+\alpha_2) \neq 0
\end{align*}
Further, \textit{conjectured} that:
\begin{conj}\label{thm.00.03}[\cite{07} Corollary 2.32]
The algebra $\mathcal{D}_4(\alpha_1, \alpha_2)$ is semisimple if and only if:
\begin{align*}
(\alpha_1-\alpha_2)(\alpha_1+\alpha_2) \neq 0
\end{align*} 
\end{conj}

Moreover, it has been calculated that the radical of $\mathcal{D}_4(\alpha, -\alpha)$ is generated by the commutator, that is, $\sigma [x_{12},x_{13}]$ for all $\sigma \in \mathbb{S}_4$, while that of $\mathcal{D}_4(\alpha, \alpha)$ is generated by:
\begin{align*}
\sigma(x_{12}x_{13}+ x_{12}x_{14} + x_{12}x_{23} + x_{13}x_{23} + x_{14}x_{12} + \alpha_1 ) \mid \sigma \in \mathbb{S}_4
\end{align*} 
where both ideals are of $552$-dimensional, and their corresponding quotient algebras are of $24$-dimension. \\ 

This paper is motivated by Conjecture \ref{thm.00.03} and mainly aim to address the representation theory of the non-semisimple deformation of $\mathcal{D}_4(\alpha_1, \alpha_2)$ as follows: \\
In Section \ref{section 01} we present some preliminaries by setting up common terminology, notation and elementary results. Section \ref{section 02} is dedicated to the study of $\mathcal{D}_n(\alpha, -\alpha)$ where we start by proving that the algebra $\mathcal{D}_n(\alpha, -\alpha)$ admits a quiver-presentation that coincidences with the graph of the nil-Coxeter group associated with $\mathbb{S}_n$, we then consider the special case of $n= 4$ where we prove that generic indecomposable projective $\mathcal{D}_4(\alpha, -\alpha)$-modules are radically graded and are isomorphic to the nil-Coxeter algebra of $\mathbb{S}_4$. Finally we consider the algebra $\mathcal{D}_4(\alpha, \alpha)$ in Section \ref{section 03}, which we show to be non-basic, implying the existence of a basic associated Morita equivalent algebra, we construct an algebraic presentation of its graded generic indecomposable projective modules and then proceed to propose a family of which the none-graded version belongs to. \\

\textbf{Acknowledgments.} The first-named author is sponsored by the German Academic Exchange Service (DAAD) to which he is immensely grateful.

\section{Preliminaries} \label{section 01}

By a $K$-algebra, we mean an associative unital algebra over $K$. We say that a $K$-algebra $A$ is connected if it is not isomorphic to a direct product of two non-trivial algebras. We also denote $A$'s Jacobson's radical, that is, the intersection of all maximal ideals of $A$ by $radA$.  \\

We further understand $Ext^1_A(\rho, \rho')$ the space of extensions between two $A$-representations $\rho, \rho'$ as equivalence classes $Z^1(\rho, \rho')/B^1(\rho, \rho')$, where: $Z^1(\rho, \rho')$ denotes the space of $(1)$-cocycles, that is:
\begin{align*}
Z^1(\rho, \rho') = \{f: A \rightarrow Hom_K(\rho, \rho') \mid f(xy)=\rho(x)f(y)+f(x)\rho'(y)\}
\end{align*}
and $B^1(\rho, \rho')$ denotes the space of coboundaries.

\subsection{Graded and filtered algebras.}

\begin{defn}
Given $A$ a filtered algebra, that is, an algebra with a family of subspaces $\{F_{i+1} \subseteq  F_i, i \geq 0 \}$ such that: $1 \in F_0$, $F_iF_j \subseteq F_{i+j}$ and $\cup F_n= A$. We define the associated graded algebra of $A$, denoted by $grA$ by setting $(grA)_n= F_n/F_{n+1}$ and $grA= \oplus (grA)_n$. 
\end{defn}

\begin{exmp}
Given $A$ a finite-dimensional $K$-algebra. Then $A$ is \textit{radically} filtered as follows:
\begin{align*}
F_r= (radA)^r \subset \cdots F_2= (radA)^2 \subset F= radA \subset F_0= A
\end{align*}
where $r$ is the minimal positive integer such that $F_{r+1}= 0$\footnote{such $r$ exists for Jacobson radical are nilpotent}
\end{exmp}

\begin{note}
Unless otherwise mentioned, $grA$ denotes the associated graded algebra of a $K$-algebra $A$ with respect to the radical filtration.
\end{note}

\begin{defn}
Given $\phi: M \rightarrow N$ a filtered homomorphism, that is, $\phi(M_j) \subseteq \phi(M)\cap N_j$. If it happens that $\phi(M_j)= \phi(M) \cap N_j$ for each $j$ applicable, then $\phi$ is called \textit{strict}.
\end{defn}

\begin{exmp}
If $\alpha: M \rightarrow N$ is an arbitrary homomorphism and $M$ is given the induced filtration $M_j= \alpha^-(\alpha(M) \cap N_j)$ then $\alpha$ is a strict filtered homomorphism. Similarly, for $\alpha$ surjective and if $N$ is given the induced filtration $N_j= \alpha(M_j)$, then $\alpha$ is strict as well. 
\end{exmp}

\begin{cor}[\cite{08} Corollary 6.14]
For $\phi: M \rightarrow N$ a filtered homomorphism. Then $gr\phi$ is injective (surjective) if and only if $\phi$ is injective (surjective) and $\phi$ is strict.
\end{cor}

\subsection{Basic algebras.}

\begin{defn}
A finite-dimensional $K$-algebra $A$ is said to be basic if and only if the quotient algebra $A/rad(A)$ is isomorphic to a product of copies of $K$.
\end{defn}

\begin{rem}
Every simple module over a basic $K$-algebra is one-dimensional.
\end{rem}

\begin{defn}[\cite{09} Corollary 6.10]
Let $A$ be a (not necessarily basic) $K$-algebra, then there exists a basic $K$-algebra $A^b$ associated with $A$ such that is \textit{Morita equivalent} to $A$, that is, there exists a $K$-linear equivalence of the modules categories $modA^b$ and $modA$.
\end{defn}

\subsection{The quiver of a finite dimensional algebra}

\begin{defn}
A \textit{quiver} is a quadruple $Q= (Q_0, Q_1, s, t)$ with $Q_0$ and $Q_1$ finite sets and two maps $s, t: Q_1 \rightarrow Q_0$. The elements of $Q_0$ and $Q_1$ are called \textit{vertices} and \textit{arrows} of $Q$ respectively. We say an arrow $\alpha$ in $Q_1$ \textit{starts} in $s(\alpha)$ and \textit{terminates} in $t(\alpha)$.
\end{defn}

\begin{exmp}
Let $A$ be a basic and connected finite-dimensional $K$-algebra and $\{S_1, \cdots, S_n \}$ a complete set of simple $A$-modules. The \textit{(ordinary) quiver} of $A$, denoted by $Q_A$, is defined as follows:
\begin{enumerate}
\item The vertices of $Q_A$ are numbers $1, \cdots, n$ which are in bijective correspondence with the simples $S_1, \cdots, S_n$.
\item Given two points $i, j \in Q_A$, the arrows $\alpha: i \rightarrow j$ are in bijective correspondence with the vectors in a basis of the $K$-vector space $Ext^1_A(S_i, S_j)$.
\end{enumerate}
\end{exmp}

\begin{rem}
A \textit{path} of length $m \geq 1$ in $Q$ is a tuple $(\alpha_1, \cdots, \alpha_m)$ of arrows of $Q$ such that $s(\alpha_i)= t(\alpha_{i+1})$ for all $1 \leq i \leq m-1$, we write such path as $\alpha_1\cdots\alpha_m$ if no misunderstanding occurs. Additionally, for each vertex $i$ of $Q$ there exists a path $e_i$ of trivial length such that $s(e_i)= t(e_i)= i$.
\end{rem}

\begin{defn}
Let $Q$ be a quiver. The \textit{path algebra} $KQ$ of $Q$ is the $K$-algebra whose underlying $K$-vector space has as its basis the set of all $Q$-paths of length $l \geq 0$ in $Q$ such that the product of two basis vectors $\alpha_1\cdots\alpha_{m}$ and $\beta_1\cdots\beta_{m'}$ is trivial if $t(\alpha_m) \neq s(\beta_1)$ and equal to the composed path $\alpha_1\cdots\alpha_m\beta_1\cdots\beta_{m'}$ otherwise.
\end{defn}

\begin{thm}[Gabriel \cite{09} Theorem 3.7]
Let $A$ be a basic and connected finite-dimensional $K$-algebra. There exists an admissible ideal $I$ of $KQ_A$ such that $A \cong KQ_A/I$.
\end{thm}

\section{Part 01: Representation theory of $\mathcal{D}_n(\alpha, -\alpha)$.} \label{section 02}

For convenience, we denote $\Lambda \coloneqq \mathcal{D}_n(\alpha, -\alpha)$, we further normalize the $K$ parameter $\alpha$ to $1_K$.

\begin{lem} \label{thm.02.01}
Given $\sigma \in \mathbb{S}_n$. The algebra homomorphism: $\rho_\sigma: A \rightarrow K$ defined by mapping a generator $x_{ij}$ - with $1 \leq i \neq j \leq n$- to:
\begin{align*}
\rho_\sigma(x_{ij}) = \begin{cases} 
+1 \mid \sigma(i) < \sigma(j), \\
-1 \mid  \sigma(i) > \sigma(j),
\end{cases}
\end{align*}
is a well defined one-dimensional representation of $\Lambda$.
\begin{proof}
The proof follows verifying that subjecting $\rho_\sigma$ to the defining relations of $\Lambda$ yields valid equations in $K$: Indeed as (\ref{eq.00.04}) and (\ref{eq.00.05}) hold trivially, one only has to check that for $1 \leq i, j, k \leq n$ distinct, then:
\begin{align*}
\rho_\sigma(x_{ij})\rho_\sigma(x_{jk})+ \rho_\sigma(x_{jk})\rho_\sigma(x_{ki}) + \rho_\sigma(x_{ki})\rho_\sigma(x_{ij}) &= -1
\end{align*}
which is easily verified as assuming that $1 \leq i,  j, k \leq n$ distinct, yield that the inequality regarding $\sigma(i), \sigma(j), \sigma(k)$ has exactly one of six possibilities:
\begin{align*}
\sigma(i) < \sigma(j) < \sigma(k). && \sigma(i) < \sigma(k) < \sigma(j). \\
\sigma(j) < \sigma(i) < \sigma(k). && \sigma(j) < \sigma(k) < \sigma(i). \\
\sigma(k) < \sigma(i) < \sigma(j). && \sigma(k) < \sigma(j) < \sigma(i)
\end{align*}
\end{proof}
\end{lem}

\begin{lem}  \label{thm.02.02}
Given $\rho$ a one-dimensional $\mathcal{D}_n(\alpha, -\alpha)$-representation. Then there exists $\sigma \in \mathbb{S}_n$ such that $\rho = \rho_\sigma$.
\begin{proof}
Let $\rho$ be a one-dimensional $\Lambda$-representation, (\ref{eq.00.04}) implies that $\rho(x_{ij})= \pm1$, furthermore, for distinct $1 \leq i, j, k \leq n$, then (\ref{eq.00.06}) can hold if one of the following possibilities occurs:
\begin{align*}
+y_{ij}= +y_{jk}= +y_{il}= +1. && +y_{ij}= -y_{jk}= +y_{ik}= +1. \\
-y_{ij}= +y_{jk}= +y_{ik}= +1. && +y_{ij}= +y_{jk}= -y_{ik}= +1. \\
+y_{ij}= -y_{jk}= +y_{ik}= +1. && -y_{ij}= -y_{jk}= -y_{ik}= +1. 
\end{align*}
Which in and of itself yield the claim.
\end{proof}
\end{lem}

\begin{rem}  \label{rem.02.01}
The previous lemma can be alternatively proven by setting: $l_i \coloneqq |\mathbb{L}_i|$, $r_i \coloneqq  |\mathbb{R}_i|$, where: 
\begin{align*}
\mathbb{L}_i &\coloneqq \{1 \leq j \leq i-1 \mid \rho(x_{ji})= -1 \} \subseteq \{0, \cdots, i-1 \} \\
\mathbb{R}_i &\coloneqq  \{i+1 \leq j \leq n \mid \rho(x_{ij})= -1 \} \subseteq  \{0, \cdots, n-i \}
\end{align*}
and defining the mapping $\sigma$ on $\{1, \cdots, n \}$ where $\sigma(i)= i+r_i-l_i$, the claim follows by showing that $\sigma$ is indeed a permutation such that $\rho= \rho_\sigma$.
\end{rem}

\begin{lem}  \label{thm.02.03}
Given $\sigma, \tau \in \mathbb{S}_n$, then for all $x \in A$ we have:
\begin{align*}
\rho_\sigma(\tau(x)) &= \rho_{\sigma\tau}(x).
\end{align*}
\begin{proof}
Since $\rho_\sigma, \rho_\tau$ and the group action of $\mathbb{S}_n$ is multiplicative, it is enough to verify the claim for a generator $x_{ij}$ with $1 \leq i \neq j \leq n$ which hold directly since:
\begin{align*}
\rho_\sigma(\tau(x_{ij})) = \rho_{\sigma}(x_{\tau(i)\tau(j)}) = \rho_{\sigma\tau}(x_{ij}).
\end{align*}
\end{proof}
\end{lem}

\begin{quote}
\textbf{Question:} Computing the radical of $\Lambda$ for a positive integer $n \geq 5$ remains open for the moment. We highlight that our results would imply the basicness of the algebra $\Lambda$ once we were able to verify that $\{\rho_\sigma \mid \sigma \in \mathbb{S}_n\}$ gives a complete system of simple $\Lambda$-representations.
\end{quote}

Our aim in at this point is prove the following theorem:
\begin{thm}  \label{thm.02.04}
Given $\sigma, \tau \in \mathbb{S}_n$. If $\tau = \overline{g}\sigma$ where $\overline{g}$ denotes a non-simple transposition of $\mathbb{S}_n$, then $dim_KExt^1_\Lambda(\rho_\sigma, \rho_\tau) = 1$, and $0$ in any other case.
\end{thm}

\begin{rem} \label{rem.02.02}
In the purpose of proving Theorem \ref{thm.02.04}, we start by utilizing Lemma \ref{thm.02.03}, which implies that we may set $\sigma= e$ with no further restrictions. Furthermore, for $\tau \in \mathbb{S}_n$, the space of extensions $Ext^1_\Lambda(\rho_e, \rho_\tau)$ has a generating set of the form: 
\begin{align*}
\{f_{ij}+B^1(\rho_e, \rho_\tau)  \mid 1 \leq i<j \leq n \mid f_{ij} \in K\}
\end{align*}
such that the following hold: 
\begin{subequations} 
\begin{align} 
f_{ij} (1 + \rho_\tau(x_{ij})) &= 0 \label{eq.00.07}  \\
f_{ij} (1 - \rho_\tau(x_{kl})) - f_{kl} (1 - \rho_\tau(x_{ij})) &= 0 \label{eq.00.08} \\
f_{ij}(\rho_\tau(x_{jk}) - 1) + f_{jk} (1 - \rho_\tau(x_{ik})) - f_{ik} (1 + \rho_\tau(x_{ij})) &= 0  \label{eq.00.09}\\
f_{ij}(1 - \rho_\tau(x_{ik})) + f_{jk} (\rho_\tau(x_{ij}) - 1) - f_{ik}(\rho_\tau(x_{jk}) + 1) &= 0 \label{eq.00.10}
\end{align}
\end{subequations}
where $1 \leq i < j \leq n$ in \ref{eq.00.07}, $1 \leq i, j, k, l \leq n$ distinct in \ref{eq.00.08}, and $1 \leq i < j < k \leq n$ in both \ref{eq.00.09} and \ref{eq.00.10}. 
\end{rem}

We now verify \ref{thm.02.04} by showcasing the following set of propositions:

\begin{prop}  \label{thm.02.05}
If $\tau$ is a non-simple transposition. Then:
\begin{align*}
dim_kExt^1_\Lambda(\rho_{e}, \rho_\tau)= 1
\end{align*}
\begin{proof}
Assuming that $\tau= (s, t)$ a non-simple transposition, then $\rho_\tau(x_{ij}) = -1$ if and only if:
\begin{align*}
(i = s \text{ and } j \leq t) \text{ or } (i \geq s \text{ and }j = t)
\end{align*}
Therefore, for all $i < j$, \ref{eq.00.07} implies that $f_{ij} = 0$ except those of the form:
\begin{align*}
f_{s(s+1)} \cdots f_{st} \text{ and } f_{(s+1)t} \cdots f_{(t-1)t}
\end{align*}
Further, we have:
\begin{align*}
f_{s(s+r)} = f_{(s+r)t} &\text{ via } \ref{eq.00.09} \text{ } (i= s, j= s+r, k= t) &\mid 1 \leq r \leq s-1. \\
f_{s(s+1)} = f_{s(s+r)} &\text{ via } \ref{eq.00.10} \text{ } (i= s, j= s+1, k= s+r) &\mid 2 \leq r \leq s-1.
\end{align*}
Therefore, we are in the situation where: 
\begin{align*}
f_{s(s+1)} = \cdots = f_{s(t-1)} = f_{(s+1)t} = \cdots = f_{(t-1)t} .
\end{align*}
Which with a proper choice of basis assert the claim.
\end{proof}
\end{prop}

\begin{rem} \label{rem.02.03}
Given $\tau \in \mathbb{S}_n$ such that $\tau \neq \overline{g}$. Then $\tau$ has one of the following form:
\begin{align*}
\begin{cases}
\tau= e \\
\tau \text{ is a simple transposition} \\
\tau \text{ is a cycle of length } p \geq 3 \\
\tau \text{ has at least two commutative cycles of length }p, q \geq 2
\end{cases}
\end{align*}
\end{rem}

\begin{prop}  \label{thm.02.06}
If $\tau= e$. Then $dim_kExt^1_\Lambda(\rho_{e}, \rho_\tau) = 0$.
\begin{proof}
Assuming that $\tau= e$, this would imply that $\rho_\tau(x_{ij})= 1$ for all $1 \leq i < j \leq n$. Now we have $f_{ij}= 0$  directly via \ref{eq.00.07} which in and of itself assures the claim.
\end{proof}
\end{prop}

\begin{prop}  \label{thm.02.07}
If $\tau$ is a simple transposition. Then $dim_kExt^1_\Lambda(\rho_{id}, \rho_\tau)$ $= 0$.
\begin{proof}
Assuming that $\tau= (s, s+1)$ for $s= 1, 2, \cdots, n$, this would imply that $\rho_\tau(x_{ij})= 1$ for all $1 \leq i < j \leq n$ except for $(i, j)= (s, s+1)$. Now we have $f_{ij}= 0$ for all $1 \leq i < j \leq n$ except for $(i, j)= (s, s+1)$ directly via \ref{eq.00.07}, which hold the claim with a proper change of basis.
\end{proof}
\end{prop}

\begin{prop}  \label{thm.02.08}
If $\tau$ is a cycle of length $p \geq 3$. Then $dim_kExt^1_\Lambda(\rho_{id}, \rho_\tau)$ $= 0$.
\begin{proof}
Assuming that $\tau= (a_1, \cdots, a_p)$ is a cycle of length $p \geq 3$ where $1 \leq a_1, \cdots, a_p \leq n$, $a_i \neq a_j$ for all $i \neq j$ ordered such that $a_1 < a_j$ for all $2 \leq j \leq p$.
\begin{align*}
\text{Since } a_1 < a_p \text{ and } \tau(a_p)= a_1 < a_2= \tau(a_1) \implies \rho_\tau(x_{a_1a_p})= -1.
\end{align*}
We make a bases change such that $f_{a_1a_p}= 0$. \\ \\
As for the case of $i \leq a_1$, we have:
\begin{align*}
\tau(i)= i < a_2 = \tau(a_1) &\implies \rho_\tau(x_{ia_1})= 1. & \implies f_{ia_1}= 0 &\text{ via } \ref{eq.00.07} \\
\tau(i)= i < a_1 = \tau(a_p) &\implies \rho_\tau(x_{a_pi})= 1. & \implies f_{a_pi}= 0 &\text{ via } \ref{eq.00.07} 
\end{align*}
As for the case of $a_p< i$, we have: 
\begin{align*}
\tau(a_p)= a_1 < i = \tau(i) &\implies \rho_\tau(x_{a_pi})= +1 & \implies f_{a_pi}= 0 &\text{ via } \ref{eq.00.07}  \\
\tau(a_p)= a_1 < a_2= \tau(a_p) &\implies \rho_\tau(x_{a_1a_p})= -1 &\implies f_{a_1i}= 0 &\text{ via } \ref{eq.00.09}
\end{align*}
Now if $a_p= a_1+ 1$ then all possible cases for $i$ has been considered and the claim follows. If not, then for all $a_1 < i < a_p$ we have $\tau(i)> \tau(a_p)$ and hence $f_{a_1i}= f_{ia_p}$ via \ref{eq.00.09}, and we have two cases to consider here as well: \\
If $\rho_\tau(x_{a_1i})= 1$ then $f_{a_1i}= 0$ via \ref{eq.00.07}. Otherwise, we are in the situation where:
\begin{align*}
a_1 < i < a_p \text{ and } \tau(a_1)> \tau(i) > \tau(a_p)
\end{align*}
which implies that $p > 3$\footnote{otherwise we get $i < a_p$ and $i > a_p$ a clear contradiction.} and hence $\tau$ is of the form:
\begin{align*}
\tau= (a_1, \cdots, r, \cdots, k, \cdots, a_p) \mid r > k > a_1 
\end{align*}
which implies that $f_{a_1i}= f_{ia_p}= 0$ via \ref{eq.00.08} for $(a_1 < k), (r < a_p)$. Therefore, with all possible cases considered for $\tau$ a cycle of length $p \geq 3$ we have $dim_kExt^1_\Lambda(\rho_{id}, \rho_\tau) = 0$ as claimed.
\end{proof}
\end{prop}

\begin{prop}  \label{thm.02.09}
If $\tau$ has at least two commutative cycles of length $p, q \geq 2$. Then $dim_kExt^1_\Lambda(\rho_{id}, \rho_\tau)$ $= 0$.
\begin{proof}
Assume that $\tau$ has commutative cycles $(a_1, \cdots, a_p), (b_1, \cdots, b_q)$ for $p, q \geq 2$ such that $a_1 < a_j$ for $2 \leq j \leq p$, 
$b_1 < b_j$ for $2 \leq j \leq q$ and $a_1 < b_1$. We start with:
\begin{align*}
&a_1 < a_p \text{ and } \tau(a_p)= a_1 < a_2= \tau(a_1) &\implies \rho_\tau(x_{a_1a_p})= -1 \\
&b_1 < b_q \text{ and } \tau(b_q)= b_1 < b_2= \tau(b_1) &\implies \rho_\tau(x_{b_1b_q})= -1
\end{align*}
For convenience of reference, we denote $a_1= s, a_p= t, b_1= s', b_q= t'$, and consider a change of basis such that $f_{st}= 0$, further:
\begin{align*}
\rho_\tau(x_{st})= -1 &\implies f_{kl}= 0 &&\mid 1 \leq k < l \leq n; s \neq t \neq k \neq l &\text{ via } \ref{eq.00.08}.  \\
\rho_\tau(x_{s't'})= -1 &\implies f_{kl}= 0 &&\mid 1 \leq k < l \leq n; s' \neq t' \neq k \neq l &\text{ via } \ref{eq.00.08}.
\end{align*}
Note that the first line implies in particular $f_{s't'}= 0$. And we are left with the following cases of $f_{ss'}, f_{st'}, f_{min(s',t')max(s', t')}$ and  $f_{min(t', t)max(t', t)}$. On one hand, for $f_{ss'}$, we have:
\begin{align*}
\text{if } \rho_\tau(x_{ss'})= +1 &\implies f_{ss'}= 0 \text{ vie } \ref{eq.00.07} \\
\text{if } \rho_\tau(x_{ss'})= -1 &\implies f_{ss'}= 0 \text{ vie } \ref{eq.00.09} & (i= s, j= s', k= t')
\end{align*}
While on the other hand, the remaining cases are processed by differentiating possible orderings of $t, s', t'$: \\
If $s < s' < t < t'$, then:
\begin{align*}
f_{s't}= 0 &\text{ via } \ref{eq.00.09} \text{ } (i= s, j= s', k= t), \\
f_{t't}= 0 &\text{ via } \ref{eq.00.10} \text{ } (i= s', j= t', k= t), \\
f_{st'}= 0 &\text{ via } \ref{eq.00.10} \text{ } (i= s, j= t', k= t).
\end{align*}
If $s< s' < t < t'$, then we have:
\begin{align*}
\tau(t)= \tau(a_p)= a_1 < b_1 = \tau(b_q)= \tau(t') \implies \rho_{\tau}(x_{tt'})= +1
\end{align*}
and then:
\begin{align*}
f_{tt'}= 0 &\text{ via } \ref{eq.00.07} \text{ } (i= t, j= t'), \\
f_{s't}= 0 &\text{ via } \ref{eq.00.09} \text{ } (i= s, j= s', k= t), \\
f_{st'}= 0 &\text{ via } \ref{eq.00.10} \text{ } (i= s, j= t, k= t').
\end{align*}
Finally, for the case of $s < t < s' < t'$, then we have $\\rho_\tau(x_{tt'})= +1$ and:
\begin{align*}
f_{tt'}= 0 &\text{ via } \ref{eq.00.07} \text{ } (i= t, j= t'), \\
f_{ts'}= 0 &\text{ via } \ref{eq.00.10} \text{ } (i= s, j= t, k= s'), \\
f_{st'}= 0 &\text{ via } \ref{eq.00.10} \text{ } (i= s, j= t, k= t').
\end{align*}
Therefore, with all possible cases considered for $\tau$ an $\mathbb{S}_n$-element such that it has at least two commutative cycles of length $p, q \geq 2$, we have $dim_kExt^1_\Lambda(\rho_{id}, \rho_\tau) = 0$ as claimed.
\end{proof}
\end{prop}

\begin{note}
Given $\sigma \in \mathbb{S}_n$, let $\tau= \overline{g}\sigma$ for $\overline{g}$ some non-simple transposition of $\mathbb{S}_n$. Then the space of extensions $Ext^1_\Lambda(\rho_\sigma, \rho_\tau)$ is one-dimensional by Theorem \ref{thm.02.04}, we denote the single generator of such space by $x(\sigma; \tau)$.
\end{note}

\begin{rem} \label{rem.02.04}
We conclude the connectedness of the algebra $\Lambda$ by Theorem \ref{thm.02.04} and the fact that the group $\mathbb{S}_n$ can be generated by non-simple transpositions. In other words, for all positive integers $n \geq 3$, the algebra $\mathcal{D}_n(\alpha, -\alpha)$ cannot be written as direct product of two non-trivial algebras.
\end{rem}

\textbf{The special case of $n= 4$.} We consider the special case of $\Lambda \coloneqq \mathcal{D}_4(\alpha, -\alpha)$ where the $K$ parameter $\alpha$ remain normalized to $1_K$.  Furthermore, we set $t_1= (1, 3)$, $t_2= (1, 4)$ and $t_3= (2, 4)$ the non-simple transpositions of $\mathbb{S}_4$.

\subsubsection{Quiver-presentation of $\Lambda$.}

As the Jacobson radical of the finite-dimensional algebra $\Lambda$ is generated by the commutator, we deduce that the algebra $\Lambda$ is basic with a complete system of simple representations: $\{\rho_\sigma \mid \sigma \in \mathbb{S}_4 \}$, moreover, Remark \ref{rem.02.04} implies that the algebra $\Lambda$ is connected. Therefore, by Theorem \ref{thm.02.04} we conclude that the ordinary quiver of $\Lambda$ denoted by $Q_\Lambda$ has a vertices set of the form $\{\sigma \mid \sigma \in \mathbb{S}_4 \}$, and there exists an arrow from $\sigma$ to $\tau$ labeled by $\alpha(\sigma; \tau)$ if and only if $\tau= t_i.\sigma$ for $i= 1, 2, 3$. \\

Denote by $\phi$ Gabriel's theorem morphism associated with $\Lambda$. By Gabriel's Theorem we deduce that $\Lambda \cong KQ_\Lambda/ker(\phi)$, furthermore, we observe that:
\begin{align*}
KQ_\Lambda&= \oplus_{\sigma \in \mathbb{S}_4} e_\sigma KQ_\Lambda
\end{align*}
This in particular implies that from the viewpoint of representation theory, the study of $\Lambda$ can be reduced to that of $\Gamma$ the indecomposable projective $\Lambda$-representation understood to be a quotient of $e_e KQ_\Lambda$ by the kernel of $\pi= \phi_{|_{e_e}}$.  

\begin{rem} \label{rem.02.05}
When considered in terms of Gabriel's theorem, the action of the symmetric group $\mathbb{S}_4$ on $\Lambda$ is understood for any two arrows $\alpha_1, \alpha_2$ as $\sigma(\alpha_1\alpha_2)= \sigma\alpha_1\sigma\alpha_2$ where:
\begin{align*}
\sigma\alpha(\tau_1; \tau_2)= \alpha(\sigma\tau_1; \sigma\tau_2) &&\mid \tau_1, \tau_2, \sigma \in \mathbb{S}_4
\end{align*}
\end{rem}

\subsubsection{Quiver-representation of $gr\Gamma$.}

\begin{prop}\label{thm.02.10}
\begin{align*}
r_i \coloneqq &\alpha(e; t_i)\alpha(t_i; e) \in ker(gr\pi) && \mid i= 1, 2, 3.
\end{align*}
\begin{proof}
Consider the $e_eKQ_\Lambda$-module $M$ given as the quotient by the (two-sided) ideal generated by:
\begin{align*}
\{e_erad^3KQ_\Lambda, rad^2KQ_\Lambda e_\sigma \mid e \neq \sigma \in \mathbb{S}_4 \}
\end{align*}
we remark here that $M$ exists as a $e_eKQ_\Lambda/ker(gr\pi)$-module if and only if the following (graded) algebra map:
\begin{align*}
\rho &: \Lambda \rightarrow End(K^{5})  \\
& x_{st} \mapsto \rho(x_{st})=  \begin{bmatrix} 
\rho_{e}(x_{st})& x_{st}(e; t_1)& x_{st}(e; t_2)& x_{st}(e; t_3)& g_{st(e;e)} \\
0&  \rho_{t_1}(x_{st})& 0& 0& x_{st}(t_1; e) \\
0& 0& \rho_{t_2}(x_{st})& 0& x_{st}(t_2; e) \\
0& 0& 0& \rho_{t_3}(x_{st})& x_{st}(t_3; e)  \\
0& 0& 0& 0& \rho_{e}(x_{st})
\end{bmatrix}
\end{align*}
exists\footnote{where $g_{st(e;e)}$ a set of $K$-parameters determined by the defining relations of $\Lambda$} \textbf{up to the third power of the radical} as a $\Lambda$-representation, that is, if and only if:
\begin{align*}
x_{ij}(e; t_1)x_{ij}(t_i; e)+rad^3\Lambda&= 0 \mid i= 1, 2, 3. \text{ per } \ref{eq.00.04}, \ref{eq.00.06}
\end{align*}
that is, if and only if: 
\begin{align*}
gr\pi(\alpha(e; t_i)\alpha(t_i; e))&=0 \mid i= 1, 2, 3. 
\end{align*}
asserting the claim.
\end{proof}
\end{prop}

An identical method of argument as before proposes the following:

\begin{prop} \label{thm.02.11}
The following is of $ker(gr\pi) $
\begin{align*}
r_3 \coloneqq &\alpha(e; t_1)\alpha(t_1; t_1t_3)-\alpha(e; t_3)\alpha(t_3; t_1t_3)  \\
r_4 \coloneqq &\alpha(e; t_1)\alpha(t_1; t_{1}t_2)\alpha(t_{1}t_2; t_{1}t_2t_1)-\alpha(e; t_2)\alpha(g_2; t_{2}t_1)\alpha(t_{2}t_1; t_{2}t_1t_2) 
\end{align*}
\end{prop}

\begin{rem} \label{rem.02.06}
Our computations suggest that there exists no further non-trivial relations of length $n \geq 4$. In other words, we conclude that:
\begin{align*}
ker(gr\pi)= \{\sigma.r_i, \mid \sigma \in \mathbb{S}_4, i= 1, \cdots, 4. \}
\end{align*}
\end{rem}

For $i= 1, 2, 3$, we rename paths $\alpha(\sigma; t_i\sigma)$ to $s_i$. This, along corresponding paths composition to the obvious multiplication implies that we may further identify the algebra $gr\Gamma$ with that of the bounded free associative algebra $K\langle s_1, s_2, s_3 \rangle$. In other words:

\begin{cor} \label{thm.02.12}
Up to a higher power of the radical, the algebra $gr\Gamma$ is isomorphic to the nil-Coxeter algebra associated with $\mathbb{S}_4$, that is:
\begin{align*}
gr\Gamma \cong K\langle s_1, s_2, s_3 \rangle/ker(gr\pi) \mid ker(gr\pi)= \begin{cases}
s^2_i &\mid i= 1, 2, 3 \\
s_1s_3-s_3s_1 \\
s_is_2s_i-s_2s_is_2  &\mid i= 1, 3
\end{cases}
\end{align*}
\end{cor}

\subsubsection{Quiver-representation of $\Gamma$.}

\begin{prop} \label{thm.02.13}
The algebra $\Gamma$ is isomorphic to the nil-Coxeter algebra associated with $\mathbb{S}_4$.
\begin{proof}
This is a direct consequence of the deletion property of Coxeter systems combined with the fact that elements of the ideal $ker(gr\pi)$ are minimal and maximal in the precise length sense.
\end{proof}
\end{prop}

\section{Part 02: Representation theory of $\mathcal{D}_4(\alpha, +\alpha)$.} \label{section 03}

For convenience, we denote $\Lambda \coloneqq \mathcal{D}_4(\alpha, \alpha)$, the $K$-parameter $\alpha$ remains normalized to $1_K$.

\begin{lem} \label{thm.03.01}
The algebra $\Lambda$ is not basic.
\begin{proof}
Assume for a contradiction that the algebra $\Lambda$ is basic, that is, every simple $\Lambda$-module is one-dimensional say of the form:
\begin{align*}
\rho &: \Lambda \rightarrow K  \\
& x_{ij} \mapsto \rho(x_{ij})= y_{ij}
\end{align*}
On one hand, \ref{eq.00.04} implies that $y^2_{ij} = 1$, that is, $y_{ij} = \pm 1$ for all $1 \leq i < j \leq 4$. On the other hand \ref{eq.00.06} implies that $y_{ij}y_{jk} - (y_{jk}y_{ik} + y_{ik}y_{ij}) = 1$ for $1 \leq i < j < k < 4$. We discuss: \\
As for the case of $y_{ij}y_{jk} = +1$ we get that:
\begin{align*}
y_{ij} = y_{jk} \text{ and } y_{jk}y_{ik} + y_{ik}y_{ij} = 0 
\end{align*}
While the case of $y_{ij}y_{jk} = -1$ implies: 
\begin{align*}
y_{ij} = -y_{jk} \text{ and } y_{jk}y_{ik} + y_{ik}y_{ij} = -2
\end{align*}
A clear contradiction on both cases proving that indeed, the algebra $\Lambda$ is no basic.
\end{proof}
\end{lem}

\begin{rem} \label{rem.03.01}
Lemma \ref{thm.03.01} implies the existence of a basic $K$-algebra $\Lambda^b$ such that Morita equivalent to $\Lambda$, that is: 
\begin{center}
\begin{tikzcd}
mod \Lambda \arrow{rr}{\mathfrak{F}} \arrow[swap]{rr}{\simeq}  &&mod\Lambda^b \arrow[swap]{rr}{\simeq}\arrow{rr}{\mathfrak{G}} && mod \Lambda
\end{tikzcd}
\end{center}
In particular, such equivalence preserves simplicity and exactness.
\end{rem}

\begin{prop} \label{thm.03.02}
Given $\sigma \in \mathbb{S}_4$. The algebra homomorphism $\rho_{\sigma}: \Lambda \rightarrow K^{2 \times 2}$ defined by mapping a generator $x_{ij}$ for $1 \leq i \neq j \leq 4$ to:
\begin{align*}
\rho_{\sigma}(x_{ij})= \begin{cases}
\begin{bmatrix} +1& 0 \\ 0& -1 \end{bmatrix} \mid \sigma(i, j) \in \{(1, 3), (4, 2) \} \\
\begin{bmatrix} 0& +1 \\ +1& 0 \end{bmatrix} \mid \sigma(i, j) \in \{(4, 1),  (s, s+1) \mid s= 1, 2, 3  \}  
\end{cases}
\end{align*}
is a well-defined two-dimensional simple $\Lambda$-representation.
\begin{proof}
For the claim to hold, one must verify that subjecting $\rho_{\sigma}$ to the defining relations does not yield any contradictions for all $\sigma \in \mathbb{S}_4$. As \ref{eq.00.04}, and \ref{eq.00.05} holds directly, one remark that \ref{eq.00.06} holds by observing that for $1 \leq \sigma(i) \neq \sigma(j) \neq \sigma(k) \leq 4$, then one of the following situations occurs:
\begin{align*}
\rho_{\sigma}(x_{ij}.x_{jl})&= 1 &&\rho_{\sigma}(x_{jk}x_{ki}+x_{ki}x_{ij})= 0 \\
\rho_{\sigma}(x_{jk}.x_{ki})&= 1  &&\rho_{\sigma}(x_{ij}x_{jk}+x_{ki}x_{ij})= 0 \\
\rho_{\sigma}(x_{ki}.x_{ij})&= 1 &&\rho_{\sigma}(x_{ij}x_{jk}+x_{jk}x_{ki})= 0
\end{align*}
which in and of itself assert the claim.
\end{proof}
\end{prop}

\begin{lem} \label{thm.03.03}
Let $\tau, \sigma \in \mathbb{S}_4$. Then for all $x \in \Lambda$ we have:
\begin{align*}
(\tau.\rho_\sigma)(x)= (\rho_\sigma(\tau^-.x))= \rho_{\sigma.\tau^-}(x).
\end{align*}
\begin{proof}
As $\rho_\tau$, $\rho_{\tau\sigma}$ and the group action by $\mathbb{S}_4$ is multiplicative, the claim is asserted by remarking that:
\begin{align*}
(\tau.\rho_\sigma)(x_{ij})= (\rho_\sigma(\tau^-.x_{ij}))= \rho_{\sigma\tau^-}(x_{ij}).
\end{align*}
for all $x_{ij}$ generating $\Lambda$. 
\end{proof}
\end{lem}

\begin{quote}
\textbf{Question:} On the existence of $\mu$ a 2-dimensional simple $\mathcal{D}_5(+1, +1)$-representation such that $\mu(x_{12})= \rho_{e}(x_{12})$. We deduce that $\mu(x_{45})= \pm\mu(x_{12})$, which in and of itself implies that $\mu(x_{23})= \pm\mu(x_{13})= \pm\mu(x_{12})$ by \ref{eq.00.04} and \ref{eq.00.05}  contradicting \ref{eq.00.06}. In other words $\{\rho_\sigma \mid \sigma \in \mathbb{S}_4 \}$ does not extents to a complete system of simple 2-dimensional $\mathcal{D}_5(+1, +1)$-representations, and the question of computing simple $\mathcal{D}_5(+1, +1)$-representations remains open for the moment.
\end{quote}

\begin{note}
Denote by $\mathbb{V}$ the Klein four-subgroup of the symmetric group $\mathbb{S}_4$, that is, the subgroup generated by the permutations $\nu_1= (13)(24)$ and $\nu_2= (14)(23)$, further, we set $\nu_3= \nu_1\nu_2$. One may assume with no loss of generality that $\mathbb{V}$ fixes 1 to realize the associated quotient group as the symmetric group on three letters $\{ 2, 3, 4\}$ which we denote by $\overline{\mathbb{S}_3}$.
\end{note}

\begin{prop} \label{thm.03.04}
Let $\sigma \in \mathbb{S}_4$, then $\rho_\sigma \cong \nu_i.\rho_\sigma$ via conjugation with $m_{\sigma\nu_i\sigma^-}$ for $i= 1, 2, 3$ where:  
\begin{align*}
m_{\nu_1} \coloneqq \begin{bmatrix} 0& +1\\ +1& 0 \end{bmatrix} && m_{\nu_2} \coloneqq \begin{bmatrix} +1& 0\\ 0& -1 \end{bmatrix} && m_{\nu_3} \coloneqq  \begin{bmatrix} 0& +1\\ -1& 0 \end{bmatrix}
\end{align*}
\begin{proof}
Given any $x \in \Lambda$, then:
\begin{align*}
\nu_i.\rho_\sigma(x) = \rho_{\sigma \nu^-_i}(x)= \rho_{\sigma \nu_i \sigma^- \sigma}(x)= \sigma^-.\rho_{\sigma \nu_i \sigma^-}(x).
\end{align*}
Now $\sigma \nu_i \sigma^- \in \mathbb{V}$ since $\mathbb{V}$ is a normal subgroup of $\mathbb{S}_4$, in other words, $\sigma \nu_i \sigma^- = \nu_j$ for some $j= 1, 2, 3$. The claim is then asserted since $\rho_{\nu_j} \cong \rho_e$ via $m_{\nu_j}$, that is:
\begin{align*}
\rho_{\nu_j}(x)=  m_{\nu_j}. \rho_e(x). m^-_{\nu_j} && \mid  j= 1, 2, 3.
\end{align*}
\end{proof}
\end{prop}

Our aim in at this point is prove the following theorem:

\begin{thm} \label{thm.03.05}
Let $\sigma, \tau \in \overline{\mathbb{S}_3}$, then:
\begin{align*}
Ext^1_\Lambda(\rho_\sigma, \rho_\tau) \cong Ext^1_\Lambda(\mu\rho_\sigma, \mu\rho_\tau) &&\mid \mu \in \mathbb{S}_4.
\end{align*}
Furthermore, we have:
\begin{align*}
dim_KExt^1_\Lambda(\rho_\sigma, \rho_\tau)= \begin{cases}
2 &\mid \tau = \sigma.s^-_2 \\
1 &\mid \tau= \sigma.s^-_3 \\
0 &\mid \text{ otherwise.}
\end{cases}
\end{align*}
\end{thm}

We discuss the first part of the theorem as follows:

\begin{prop} \label{thm.03.06}
Let $\sigma, \tau \in \overline{\mathbb{S}_3}$. Then:
\begin{align*}
Ext^1_\Lambda(\rho_\sigma, \rho_\tau) \cong Ext^1_\Lambda(\mu\rho_\sigma, \mu\rho_\tau) &&\mid\mu \in \mathbb{S}_4.
\end{align*}
\begin{proof}
The claim is a natural consequence of Lemma \ref{thm.03.03} induced by the group action. In details, given $f \in Z^1(\rho_\sigma, \rho_\tau)$ and $\mu \in \mathbb{S}_4$ then for all $x_{ij}$ generating $\Lambda$ we have:
\begin{align*}
\mu(f) \in Z^1(\mu\rho_\sigma, \mu\rho_\tau) \mid \mu f(x_{ij})= f(\mu^-x_{ij}) 
\end{align*}
a well-defined (1-)cocycles, that is a coboundry if and only $f$ is.
\end{proof}
\end{prop}

\begin{rem} \label{rem.03.02}
The special case of $\mu \in \mathbb{V}$ in Proposition \ref{thm.03.06} implies that the isomorphism as described is a self-inverse. Furthermore, such isomorphism can be realized by means of Proposition \ref{thm.03.04}, which implies that for an arbitrary $f \in Z^1(\rho_\sigma, \rho_\tau)$, then for all $x_{ij}$ generating $\Lambda$ we have:
\begin{align*}
v(f) \in Z^1(\mu\rho_\sigma, \mu\rho_\tau) \mid v(f)(x_{ij})= m^-_{\sigma\mu\sigma^-}f(x_{ij})m_{\tau\mu\tau^-} 
\end{align*}
a well-defined (1-)cocycles that is a coboundry if and only $f$ is.
\end{rem}

\begin{note}
We understand arbitrary $\overline{f}\in Ext^1_\Lambda(\rho_{e}, \sigma\rho_{e})$, as a class of the form $f+B^1(\rho_e,\sigma\rho_e)$ where $f: \Lambda \rightarrow K^{2.2}$ maps a generator $x_{ij}$ to:
\begin{align*}
f(x_{ij})=\begin{bmatrix} a_{ij}& b_{ij} \\c_{ij}& d_{ij} \end{bmatrix} &\mid f(xy)= \rho_e(x)f(y)+f(x)\sigma\rho_e(y).
\end{align*}
Furthermore, we denote basis changing matrices of $K^{2.2}$ by:
\begin{align*}
\lambda= \begin{bmatrix} \lambda_1& \lambda_3 \\ \lambda_2& \lambda_4 \end{bmatrix} &\mid \lambda_1, \cdots, \lambda_4 \in K.
\end{align*}
\end{note}

\begin{prop} \label{thm.03.07}
Given $\sigma \in \overline{\mathbb{S}_3}$ a non-Coxeter generator. Then the space of extensions $Ext^1_\Lambda(\rho_{e}, \sigma\rho_{e})$ is of null-dimension.
\begin{proof}
Given an arbitrary $\overline{f}\in Ext^1_\Lambda(\rho_{e}, \sigma\rho_{e})$: \\
The case of $\sigma= e$ is then resolved by a change of basis of which we set:
\begin{align*}
2\lambda_1&= 0 && 2\lambda_2= +c_{13} && 2\lambda_3= +c_{13}-2a_{12} && 2\lambda_4= +2c_{12}
\end{align*}
which would imply that:
\begin{align*}
f(x_{12})= f(x_{34})&= 0 &&\mid \text{ via } \ref{eq.00.04}, \ref{eq.00.05}, \\
f(x_{13})= f(x_{24})&= 0  &&\mid \text{ via } \ref{eq.00.04}, \ref{eq.00.06}: (i= 1, j= 2, k=3), \\
f(x_{14})= f(x_{23})&= 0 &&\mid \text{ via } \ref{eq.00.04}, \ref{eq.00.05}.
\end{align*}
The case of  $\sigma= s_2s_3s_2$ is resolved by setting:
\begin{align*}
2\lambda_1&= +a_{24} && 2\lambda_2= +c_{13} && 2\lambda_3= -b_{13} && 2\lambda_4= -d_{24}
\end{align*}
which would then imply that:
\begin{align*}
f(x_{13})= f(x_{24})&= 0 &&\mid \text{ via } \ref{eq.00.04}, \\
f(x_{12})= f(x_{23})&= 0  &&\mid \text{ via } \ref{eq.00.04}, \ref{eq.00.06}: (i= 1, j= 2, k= 3, 4), \\
f(x_{14})= f(x_{34})&=0 && \mid \text{ via } \ref{eq.00.04}, \ref{eq.00.05}.
\end{align*}
The case of $\sigma=s_3s_2$ is resolved by setting:
\begin{align*}
2\lambda_1&= +a_{14}+a_{23} && 2\lambda_2= +a_{14}-a_{23} \\
2\lambda_3&= (+a_{14}+a_{23})+2a_{13} && 2\lambda_4= (-a_{14}+a_{23})+2c_{13}
\end{align*}
which would then imply that:
\begin{align*}
f(x_{13})= f(x_{24})&= 0 &&\mid \text{ via } \ref{eq.00.04}, \\
f(x_{14})= f(x_{23})&= 0  &&\mid \text{ via } \ref{eq.00.04}, \ref{eq.00.06}: (i= 1, j= 2, k= 3, 4), \\
f(x_{34})&=0 && \mid \text{ via } \ref{eq.00.04}, \ref{eq.00.06}: (i= 1, 2, j= 3, k= 4), \\
f(x_{12})&= 0 && \mid \text{ via } \ref{eq.00.04}, \ref{eq.00.06}: (i= 1, j= 4, k= 2).
\end{align*}
Finally, the case of $\sigma= s_2s_3$ is resolved by setting:
\begin{align*}
2\lambda_1&= -a_{13}+a_{24} && 2\lambda_2= +c_{13}-c_{24} \\
2\lambda_3&= -a_{13}-a_{24} && 2\lambda_4= -c_{13}-c_{24}
\end{align*}
which would imply:
\begin{align*}
f(x_{13})= f(x_{24})&= 0 &&\mid \text{ via } \ref{eq.00.04}, \\
f(x_{12})&= 0  &&\mid \text{ via } \ref{eq.00.04}, \ref{eq.00.06}: (i= 1, j= 2, 4, k= 3, 2),\\
f(x_{34})&= 0 &&\mid \text{ via } \ref{eq.00.04}, \ref{eq.00.05}, \\
f(x_{14})&= 0  &&\mid \text{ via } \ref{eq.00.04}, \ref{eq.00.06}: (i= 1, j= 2, 4, k= 4, 3), \\
f(x_{23})&= 0  &&\mid \text{ via } \ref{eq.00.04}, \ref{eq.00.06}: (i= 2, 1, j= 3, k= 4, 2).
\end{align*}
In other words, up to an isomorphism, generic elements of the space $Ext^1_\Lambda(\rho_{e}, \rho_{\sigma})$ are trivial and therefore, the space in and of itself is of null-dimension as claimed for all $\sigma$ non-Coxeter generator.
\end{proof}
\end{prop}

\begin{prop} \label{thm.03.08}
The space of extensions $Ext^1_\Lambda(\rho_{e}, \sigma\rho_{e})$ is of one-dimensional for $\sigma= s_3$, while is two-dimensional for $\tau= s_2$.
\begin{proof}
Given an arbitrary $\overline{f}\in Ext^1_\Lambda(\rho_{e}, \sigma\rho_{e})$: \\
The case of $\sigma= s_3$ is resolved by a change of basis of which we configure:
\begin{align*}
\lambda_2= -\lambda_1+a_{14} && \lambda_3= -\lambda_1-a_{13} && \lambda_4= -\lambda_1+a_{14}-c_{13}
\end{align*}
where for $2\beta_2= -b_{12}-a_{13}+c_{13}+a_{34}$, we set:
\begin{align*}
2\lambda_1&= +b_{12}+a_{14}-c_{13}+\beta_2= -a_{13}+a_{14}+a_{34}-\beta_2
\end{align*}
This would imply that:
\begin{align*}
f(x_{13})= f(x_{24})&= 0 && \mid \text{ via } \ref{eq.00.04}, \ref{eq.00.05}, \\
f(x_{14})= f(x_{23})&= 0 && \mid \text{ via } \ref{eq.00.04}, \ref{eq.00.05}, \ref{eq.00.06}:  (i= 1, j= 2, k= 3). 
\end{align*}
Furthermore. \ref{eq.00.04} conclude that:
\begin{align*}
f(x_{12})= \beta_2.\begin{bmatrix} 0& +1 \\ -1& 0 \end{bmatrix} && e(x_{34})= \beta_2.\begin{bmatrix} +1& 0 \\ 0& +1 \end{bmatrix}
\end{align*}
Where all the other defining relations of $\Lambda$ are satisfied. In other words, up to an isomorphism, the space $Ext^1_\Lambda(\rho_{e}, s_3\rho_{e})$ has a generating set of the form: $\{\overline{f_2}= f_2+B^1(\rho_e, s_3\rho_e)\} $ where $f_2$ is defined by mapping the generators of $\Lambda$ as follows:
\begin{align*}
f_2(x_{12})&= \begin{bmatrix} 0& +1 \\ -1& 0 \end{bmatrix},
f_2(x_{34})= \begin{bmatrix} +1& 0 \\ 0& +1 \end{bmatrix},
f_2(x_{ij})= 0 \mid \text{ otherwise.}
\end{align*}
As for the case of $\sigma= s_2$, we set:
\begin{align*}
2\lambda_1&= +a_{12}-a_{34} && 2\lambda_2= -a_{12}-a_{34} \\
2\lambda_3&= -b_{12}+b_{34} && 2\lambda_4= -b_{12}-b_{34} 
\end{align*}
which would imply that:
\begin{align*}
f(x_{12})= f(x_{34})&= 0 && \mid \text{ via } \ref{eq.00.04},  \\
f(x_{13})&= 0 && \mid \text{ via } \ref{eq.00.04}, \ref{eq.00.06}: (i= 1, j= 2, 3, k= 3, 4), \\
f(x_{24})&= 0 && \mid \text{ via } \ref{eq.00.04}, \ref{eq.00.06}: (i= 1, 2, j= 4, k= 2),
\end{align*}
further:
\begin{align*}
f(x_{14})&= \beta_1.diag(+1, -1) && \mid \text{ via } \ref{eq.00.04}, \ref{eq.00.06}: (i= 1, j= 2, k= 4), \\
f(x_{23})&= \beta_3.diag(+1, +1) && \mid \text{ via } \ref{eq.00.04}, \ref{eq.00.06}: (i= 2, j= 3, k= 4).
\end{align*}
where the $K$-parameters $\beta_1, \beta_3$ are generically given as:
\begin{align*}
2\beta_1= a_{12}-b_{12}+a_{34}+b_{34}+2a_{14}, && 2\beta_3= -a_{12}-b_{12}-a_{34}+b_{34}+2a_{23}
\end{align*}
and all the other defining relations of $\Lambda$ are satisfied. In other words, up to an isomorphism, the space $Ext^1_\Lambda(\rho_{e}, s_2\rho_{e})$ has a generating set of the form: $\{\overline{f_t}= f_t+B^1(\rho_e, s_2\rho_e) \mid t= 1, 3.\} $ where $f_t$ are defined by mapping the generators of $\Lambda$ as follows:
\begin{align*}
\begin{cases}
f_1(x_{14})&= diag(+1, -1) \\
f_3(x_{23})&= diag(+1, +1)\\
f_t(x_{ij})&= 0 \mid \text{ otherwise.}
\end{cases}
\end{align*}
\end{proof}
\end{prop}

\begin{cor} \label{thm.03.08}
Let $\sigma, \tau \in \overline{\mathbb{S}_3}$. Then:
\begin{align*}
dim_KExt^1_\Lambda(\rho_\sigma, \rho_\tau)= \begin{cases}
2 &\mid \tau= \sigma.s^-_2 \\
1 &\mid \tau= \sigma.s^-_3 \\
0 &\mid \text{ otherwise.}
\end{cases}
\end{align*}
\end{cor}

\begin{rem}\label{rem.03.03}
We set $\{\overline{g_t}= g_t+B^1(\nu_1\rho_e, \nu_1\rho_{s_2}) \mid t= 1, 3.\}$ a generating set of the space $Ext^1_\Lambda(\nu_1\rho_{e}, \nu_1\rho_{s_2})$, where $g_t$ are defined by mapping the generators of $\Lambda$ as follows:
\begin{align*}
\begin{cases}
g_1(x_{14})&= diag(+1, +1) \\
g_3(x_{23})&= diag(+1, -1)\\
g_t(x_{ij})&= 0 \mid \text{ otherwise.}
\end{cases}
\end{align*}
and remark that:
\begin{align*}
\overline{f_1} \cong 
\begin{cases}
-\overline{g_3} &\text{ via } \nu_1 \\
-\overline{g_1} &\text{ via } v \mid \mu= \nu_1
\end{cases} &&
\overline{f_3} \cong \begin{cases}
-\overline{g_1} &\text{ via }\nu_1 \\
+\overline{g_3} &\text{ via } v \mid \mu= \nu_1
\end{cases}
\end{align*}
This, in particular, indicates that the induced action of the $\mathbb{V}$ (and that of $\mathbb{S}_4$ in general) does \textbf{not} extend to an action onto $\Lambda$-extensions. In other words, the algebra $\Lambda^b$ is \textbf{not} invariant under the action of $\mathbb{V}$.
\end{rem}

\begin{note}
The space of extensions $Ext^1_{\Lambda^b}(\mathfrak{F}(\rho_e), \mathfrak{F}(\rho_{s_2}))$ is two-dimensional, we set $\{\mathfrak{F}(\overline{f_1}), \mathfrak{F}(\overline{f_3})\}$ a generating set. Similarly, we set $\{\mathfrak{F}(\overline{f_2})\}$ a generating set of the one-dimensional space  $Ext^1_{\Lambda^b}(\mathfrak{F}(\rho_e), \mathfrak{F}(\rho_{s_3}))$.
\end{note}

\subsubsection{Quiver-presentation of $\Lambda^b$.}

Now that we have verified that the algebra $\Lambda^b$ is basic, connected and finite-dimensional, then Theorem \ref{thm.03.05} supported by definition implies that $Q_\Lambda$ the ordinary quiver of $\Lambda$ (and that of $\Lambda^b$ since $\Lambda$ and $\Lambda^b$ are Morita equivalent) has the following shape:
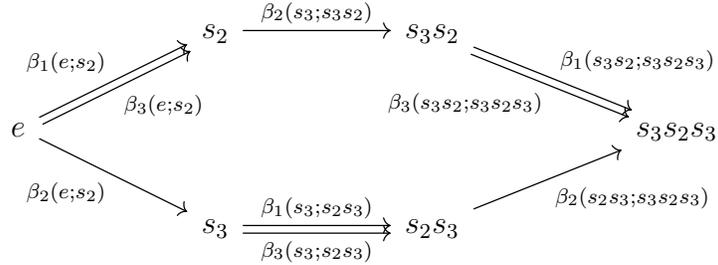
\begin{center}
\begin{tikzcd}
&& s_2 \arrow{rr}{\beta_2(s_3; s_3s_2)}&& s_3s_2 \arrow{drr}{\beta_1(s_3s_2; s_3s_2s_3)}\arrow[shift right, swap]{drr}{\beta_3(s_3s_2; s_3s_2s_3)} \\
e \arrow[swap]{drr}{\beta_2(e; s_2)} \arrow{urr}{\beta_1(e; s_2)} \arrow[shift right, swap]{urr}{\beta_3(e; s_2)} &&&&&& s_3s_2s_3 \\
&& s_3 \arrow{rr}{\beta_1(s_3; s_2s_3)}\arrow[shift right, swap]{rr}{\beta_3(s_3; s_2s_3)} && s_2s_3 \arrow[swap]{urr}{\beta_2(s_2s_3; s_3s_2s_3)}
\end{tikzcd}
\end{center}

\begin{note}
We simply write $\beta_1, \beta_2, \beta_3$ if no confusions occurs.
\end{note}

Denote by $\phi$ Gabriel's theorem morphism associated with $\Lambda^b$. By Gabriel's Theorem we deduce that $\Lambda \cong KQ_\Lambda/ker(\phi)$, furthermore, we observe that:
\begin{align*}
KQ_{\Lambda^b}&= \oplus_{\sigma \in \overline{\mathbb{S}_3}} e_\sigma KQ_\Lambda
\end{align*}
This in particular implies that from the viewpoint of representation theory, the study of $\Lambda^b$ can be reduced to that of $\Gamma$ the indecomposable projective $\Lambda^b$-representation understood to be a quotient of $e_e KQ_\Lambda$ by the kernel of $\pi= \phi_{|_{e_e}}$. 

\subsubsection{Quiver-representation of $gr\Gamma$.}

\begin{prop} \label{thm.03.09}
$gr(\pi)(r_t)= 0 \mid t= 1 \cdots, 4$, where:
\begin{align*}
r_{t}&= \beta^2_t &&\mid t= 1, 2, 3. \\
r_{4}&\coloneqq \beta_1\beta_3-\beta_3\beta_1.
\end{align*}
\begin{proof}
Consider the $e_eKQ_\Lambda$-module given as the quotient by the (two-sided) ideal generated by:
\begin{align*}
\{e_erad^3KQ_\Lambda, e_erad^2KQ_\Lambda e_\sigma \mid e \neq \sigma \in \overline{\mathbb{S}_3} \}
\end{align*}
The claim follows by verify that such module exists as an $e_eKQ_\Lambda/ker(gr\pi)$-module if and only if up to the third power of the radical we have:
\begin{align*}
\mathfrak{F}(\overline{f_i}^2)&= 0  &\mid i= 1, 2, 3. \\
\mathfrak{F}(\overline{f_1}\overline{f_3}-\overline{f_3}\overline{f_1})&= 0
\end{align*}
that is, if and only if up to the third power of the radical, we have:
\begin{align*}
\overline{f_i}^2&= 0  &\mid i= 1, 2, 3. \\
\overline{f_1}\overline{f_3}-\overline{f_3}\overline{f_1}&= 0
\end{align*}
that is, if and only if:
\begin{align*}
gr\pi(r_i)&= 0 \mid i= 1, \cdots, 4. 
\end{align*}
which asserts the claim.
\end{proof}
\end{prop}

An identical method of argument as before proposes the following:

\begin{prop} \label{thm.03.10}
$gr(\phi) r_t= 0 \mid t= 5, 6$, where:
\begin{align*}
r_5 &\coloneqq 
\beta_2\beta_1\beta_2-\beta_1\beta_2\beta_3-\beta_3\beta_2\beta_1 \\
r_6 &\coloneqq 
\beta_2\beta_3\beta_2-\beta_1\beta_2\beta_1+\beta_3\beta_2\beta_3
\end{align*}
\end{prop}

\begin{rem} \label{rem.03.04}
Our computations suggest that there exists no further non-trivial relations of length $n \geq 4$. In other words, we conclude that:
\begin{align*}
ker(gr\pi)= \{\sigma.r_t \mid \sigma \in \overline{\mathbb{S}_3}, t= 1, \cdots, 6. \}
\end{align*}
\end{rem}

identifying the path algebra $e_eKQ_\Lambda$ with the free associative algebra $K\langle s_1, s_2, s_3 \rangle$ by corresponding path composition to the obvious multiplication yield the following:

\begin{cor} \label{thm.03.12}
Up to a higher power of the radical, following hold:
\begin{align*}
gr \Gamma \cong K\langle s_1, s_2, s_3 \rangle/ker(gr(\pi))= \begin{cases}
\iota_i \coloneqq s^2_t,  \mid t= 1, 2, 3.\\
\iota_4 \coloneqq s_1s_3-s_3s_1, \\
\iota_5 \coloneqq s_2s_1s_2-s_1s_2s_3-s_3s_2s_1, \\
\iota_6 \coloneqq s_2s_3s_2-s_1s_2s_1+s_3s_2s_3.
\end{cases}
\end{align*}
\end{cor}

In particular, we find that the algebra $gr\Gamma$ is 24-dimensional. The following corollary provides a basis:

\begin{cor} \label{thm.03.13}
The following set of polynomials form a basis of the algebra $gr\Gamma$:
\begin{align*}
l&= 1: \begin{pmatrix} s_3, s_1, s_2 \end{pmatrix} &&
l= 2: \begin{pmatrix} s_3s_1,& s_3s_2,& s_1s_2,& s_2s_3,& s_2s_1 \end{pmatrix} \\
l&= 3: \begin{pmatrix} s_3s_1s_2, &s_3s_2s_3, &s_3s_2s_1, \\ s_1s_2s_3, &s_1s_2s_1, &s_2s_3s_1,  \end{pmatrix} && 
l= 4: \begin{pmatrix} s_3s_1s_2s_3, &s_3s_1s_2s_1,&  \\ s_3s_2s_3s_1, &s_1s_2s_3s_1, &s_2s_3s_1s_2,  \end{pmatrix} \\
l&= 5: \begin{pmatrix} s_3s_1s_2s_3s_1,& \\ s_3s_2s_3s_1s_2, &s_1s_2s_3s_1s_2 \end{pmatrix}
&& l= 6: \begin{pmatrix} s_3s_1s_2s_3s_1s_2  \end{pmatrix}
\end{align*}
\end{cor}

\begin{rem} \label{rem.03.05}
Remark that the basis elements $u_1 \coloneqq s_2s_3s_1s_2, u_2 \coloneqq (s_3s_1s_2)^2$ corresponds to paths $p_1, p_2$ (respectively) such that:
\begin{align*}
s(p_1)= s(p_2)= t(p_1)= t(p_2) && l(p_1), l(p_2) > 2
\end{align*}
\end{rem}

\subsubsection{Quiver-representation of $\Gamma$.}

Remark \ref{rem.03.05} along with the fact that $\Gamma \cong K\langle s_1, s_2, s_3 \rangle/ker(\pi)$ implies the existence of some $K$-polynomials $q_i$ such that:
\begin{align*}
\Gamma = K\langle s_1, s_2, s_3 \rangle/ker(\pi) = \begin{cases}
\iota_1+q_1u_1+q_{2}u_2, \\
\iota_2+q_3u_1+q_{4}u_2, \\
\iota_3+q_5u_1+q_{6}u_2, \\
\iota_4+q_7u_1+q_{8}u_2, \\
\iota_5= s_2s_1s_2-s_1s_2s_3-s_3s_2s_1, \\
\iota_6= s_2s_3s_2-s_1s_2s_1+s_3s_2s_3
\end{cases}
\end{align*}

\begin{rem}
Corollary \ref{thm.03.13} implies on one hand that:
\begin{align*}
\pi(s_1\iota_4-\iota_1s_3)&= q_7(s_1s_2s_3s_1s_2)+q_1(s_3s_2s_3s_1s_2) = 0 \\
\pi(s_1\iota_3-\iota_4s_3)&= q_5(s_1s_2s_3s_1s_2)+q_7(s_3s_2s_3s_1s_2) = 0 
\end{align*}
that is, $q_1= q_5= q_7= 0$. Further, we remark for $q$ any non-trivial $K$-polynomials that:
\begin{align*}
t_1 &\coloneqq s_1+qs_3s_2s_3s_1s_2 &&\implies t^2_1= s^2_1+2qu_2 \\
t_2 &\coloneqq s_2+qs_3s_1s_2 &&\implies t^2_2= s^2_2+qu_1+q^2u_2
\end{align*}
We remark that there exists no basis elements $u$ such that setting:
\begin{align*}
t_3\coloneqq s_3+qu &&\implies t^2_3= s^2_3+q'u_2
\end{align*}
Therefore, we set $t_3 \coloneqq s_3$ and propose:
\begin{prop}
The following hold:
\begin{align*}
\Gamma = K\langle t_1, t_2, t_3 \rangle/ker(\pi) = \begin{cases}
t^2_1, t^2_2,  t^2_3+q_1u_2, \\
t_1t_3-t_3t_1+q_2u_2, \\
t_2t_1t_2-t_1t_2t_3-t_3t_2t_1, \\
t_2t_3t_2-t_1t_2t_1+t_3t_2t_3
\end{cases}
\end{align*}
for $q_1, q_2$ two $K$-polynomials.
\end{prop}
\end{rem}

\end{document}